\documentclass[12pt]{amsart}
\usepackage{amssymb,amsfonts,amsmath,amsopn,amstext,amscd,color,
graphicx,latexsym,mathabx,mathrsfs,skak,stackrel,todonotes,xy,url,
verbatim,tikz-cd,enumerate,float}

\usepackage{hyperref}

\input xy
\xyoption{all}

\theoremstyle{remark}

\newtheorem{para}{\bf}[subsection]
\newtheorem{example}[para]{\bf Example}

\newtheorem{rem}[para]{\bf Remark}
\theoremstyle{definition}

\newtheorem{dfn}[para]{Definition}

\theoremstyle{plain}

\newtheorem{thm}[para]{Theorem}

\makeatletter
\newcommand*\bdot{{\mathpalette\bdot@{.9}}}
\newcommand*\bdot@[2]{\mathbin{\vcenter{\hbox{\scalebox{#2}{$\m@th#1\bullet$}}}}}
\makeatother

\usepackage{cleveref}
\usepackage{pdfpages}

\begin{document}

\title{Elliptic Exceptional Belyi Coverings}
\author{Cemile Kurkoglu}
\address{Denison University, OH, U.S.A.}
\email{cemile.kurkoglu@gmail.com}

\begin{abstract} An elliptic exceptional Belyi covering is a connected Belyi covering uniquely determined by its ramification scheme or the respective dessin d’enfant when the underlying compact Riemann surface has genus 1. We give our Maple algorithm and table that display our calculations for elliptic exceptional Belyi coverings up to degree 12.
\end{abstract}

\maketitle

\tableofcontents
\section{Introduction}
Let $S_g$ be a compact Riemann surface with genus $g$ and let $\mathbb{P}^1$ denote the Riemann sphere. A {\it Belyi covering} $\beta: S_g \to \mathbb{P}^1$ is a ramified covering of $\mathbb{P}^1$ with ramification points in $\{0,\,1,\,\infty\}\,.$ If $g=0$, we call $\beta$ a rational Belyi covering, if $g=1$, we call $\beta$ an elliptic Belyi covering and if $g>1$, we call $\beta$ a hyperbolic Belyi covering. If the monodromy group is $\left\langle g_{0}, g_{1}, g_{\infty}\right| g_{1} g_{2} g_{\infty}=$ $i d\rangle$, then $\lambda_{0},\,$ $\lambda_{1}\,$ and $\lambda_{\infty},$ as the cycle structures of $g_i$, are called the \textit{ramification indices} and $\left[\lambda_{0}\right]\left[\lambda_{1}\right]\left[\lambda_{\infty}\right]$ is called the \textit{ramification scheme} of $\beta\,.$\\

In Chapter 2, we revisit \cite{Ku} to give some preliminaries about Belyi coverings. We state Belyi's theorem in \cite{Belyi}. We introduce the bipartite graphs ``dessin d'enfants" with a topological structure which were studied by Grothendieck in \cite{Groth}. As proven in \cite{Zvo}, dessins and Belyi coverings are in 1-1 correspondence. We give two examples for rational Belyi coverings. Our first elliptic Belyi covering example comes from a rational one. Then we define exceptional Belyi coverings. We state a theorem that counts these coverings using the Eisenstein number. The proof of this theorem, as given in \cite{Kurk}, is based on the the works of Tutte in \cite{Tut} and the Burnside theorem. \\

Next, in Chapter 3, we illustrate with an example how it is a challenging thing to draw dessins of elliptic
exceptional Belyi coverings. We give one more example to introduce the term ``double-covering".\\

We conclude the paper with our calculations for elliptic exceptional Belyi coverings up to degree 12. We attach the Maple algorithm and the table at the end. 
\section{Exceptional Belyi Coverings}
\subsection{Belyi coverings and dessin d'enfants}
Belyi showed that $S_g$ can be defined over a number field if and only if there is on it a meromorphic function with three critical values.
\begin{thm}{\bf (Belyi)}
The following are equivalent:
\begin{enumerate}[(i)]
    \item $S_g$ is defined over $\overline{\mathbb{Q}}$.
    \item $S_g$ admits a meromorphic function $\beta: S_g \rightarrow \mathbb{P}^{1}$ with at most three ramification points.  
\end{enumerate}
\end{thm}

We also have a geometric way to represent Belyi coverings.

\begin{dfn}
 A \textit{dessin d'enfant} is a pair $(\mathcal{D}, X)$ where $X$ is an oriented compact topological surface, and $\mathcal{D} \subset X$ is a finite graph satisfying the following:

\begin{enumerate}[(i)]
    \item $\mathcal{D}$ is connected.
    \item $\mathcal{D}$ is bicoloured: two ends of every edge is colored by black and white.
    \item At each vertex is given a cyclic ordering of the edges meeting it.
    \item  $X \backslash \mathcal{D}$ is the union of finitely many topological discs, which we call \textit{faces} of $\mathcal{D}$.   
\end{enumerate}
\end{dfn}

 \begin{thm} {\bf (Zvonkin)}
    There is a one-to-one correspondence between the following:
\begin{enumerate}[(i)]
    \item $\operatorname{dessins}\left(\mathcal{D}, S_{g}\right)$
    \item Belyi coverings of $\beta: S_{g} \rightarrow \mathbb{P}^{1}$ with degree $n$
    \item the solutions of the monodromy group relation $g_{0} g_{1} g_{\infty}=i d$, where $g_{i} \in S_{n}$
\end{enumerate}
\end{thm}
The black nodes are the inverse images of 0, the white ones the inverse images of 1, and the edges of the dessin are the components of the inverse image of the line segment (0,1). The cyclic order arises from the local monodromy around the vertices.

\subsection{Examples}

\begin{example}\label{Example 2.1.6}
{\bf(Rational)} 
The power function $z\to z^n$ and the Chebyshev polynomial $T_n(\cos\varphi)=\cos n\varphi$ are two rational Belyi coverings. As classified in \cite{A09a}, the dessin d'enfants are in the form of trees. The following figure shows the corresponding dessins and the ramification schemes for these two Belyi coverings.
\newpage
\begin{figure}[htp]
\includegraphics[width=100mm]{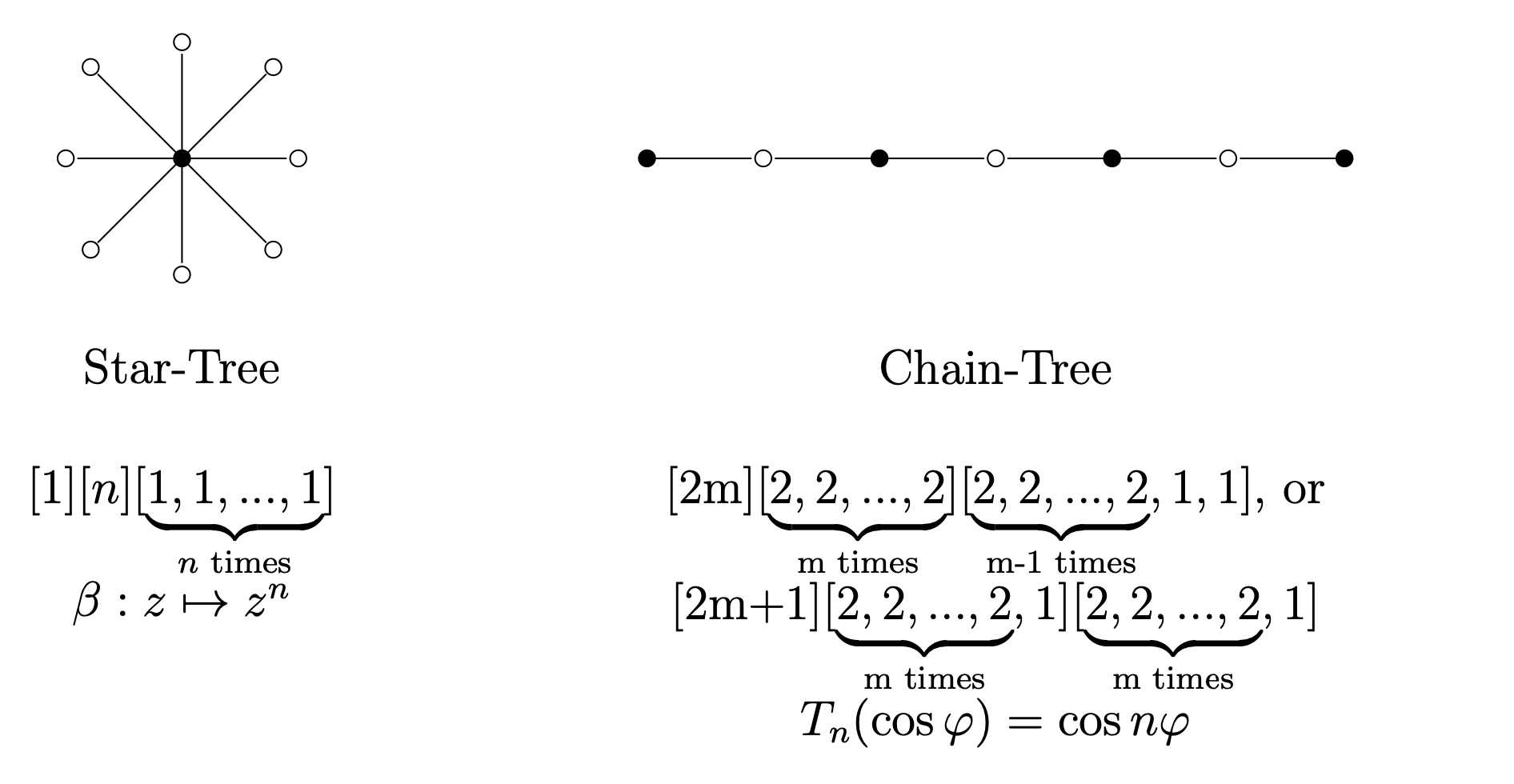}
\end{figure}
\end{example}

\begin{example}{\bf(Rational)}\label{Example 2.1.7} The following dessin is also for a rational Belyi covering but this time it is not a tree.

\begin{figure}[htp]
\includegraphics[width=80mm]{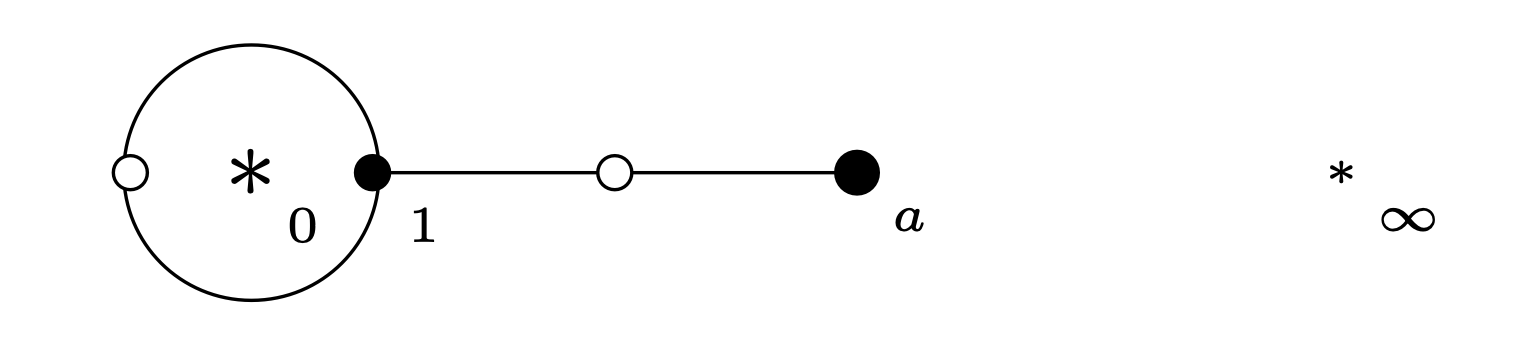}
\end{figure}

  The poles are denoted by * and as it can be seen in the dessin above, one of the poles is at the center of one of the faces and the other pole is at $\infty$. If 1 and $a$ are the roots of the corresponding function, the Belyi covering will be of the form

$$
\beta(x)=K \frac{(x-1)^{3}(x-a)}{x}
$$

since the multiplicities of the black vertices are 1 and 3. The multiplicity of the pole at $\infty$ is 3. Similarly, considering the multiplicities of the white vertices being 2,
$$
\beta(x)-1=K \frac{(x^2+bx+c)^2}{x}
$$
We find $K=-\frac{1}{64}$ and $a=9$. The roots of $\beta(x)-1$ are $3 \pm 2 \sqrt{3}$.
\end{example}

\begin{example}{\bf(Elliptic)}\label{Example 2.1.8} Now we introduce an example for the genus 1 case.
Consider the elliptic curve $$E_{1}: y^{2}=x(x-1)(x-(3+2 \sqrt{3})\,.$$ The projection on the first coordinate

\begin{align*}
\rho: E_{1} & \rightarrow \mathbb{P}^{1} \\
(x, y) & \mapsto x
\end{align*}

is a meromorphic function on $E_1$ but it is not a Belyi covering since the ramification points are $0,1,3+2 \sqrt{3}$, and $\infty$. Using the rational covering in {\bf Example \ref{Example 2.1.7}}, the composition 
$$
(x, y) \rightarrow x \rightarrow-\frac{(x-1)^{3}(x-9)}{64 x}
$$
will give us a Belyi covering. Note that $E_{1}$ has also a conjugate curve $$E_{2}: y^{2}=x(x-1)(x-(3-2 \sqrt{3})\,.$$ 
The following dessins are for $E_1$ and $E_2$ respectively.

\begin{figure}[htp]
\includegraphics[width=60mm]{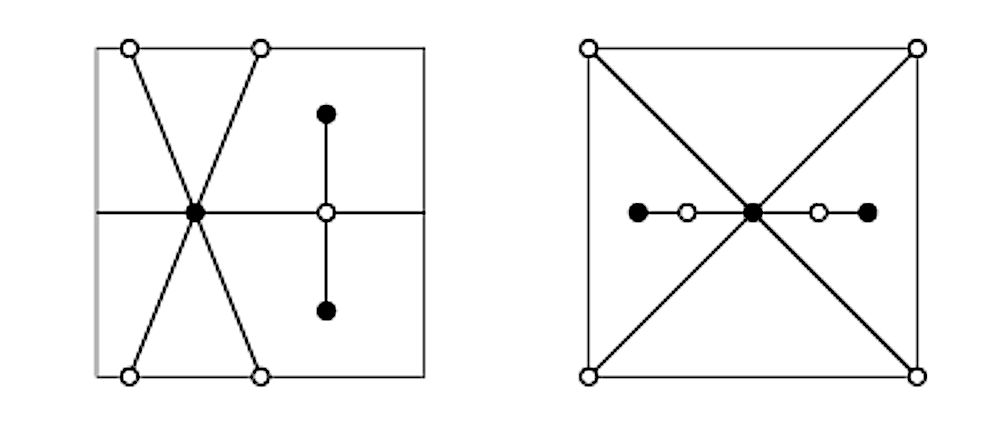}
\end{figure}

These elliptic Belyi coverings both have the same ramification scheme $[6,2][6,1,1][4,2,2]$ and they are, as algebraic curves, defined over the field $\mathbb{Q}(\sqrt{3})\,.$ 
\end{example}

\subsection{Counting exceptional Belyi coverings}
\begin{dfn}
    Let

$$
\beta: S_{g} \rightarrow \mathbb{P}^{1}
$$

be a Belyi covering. We call $\beta$ an \textit{exceptional Belyi covering} if it is uniquely determined by its ramification scheme. The corresponding dessin is also unique. 
\end{dfn}

We count the Belyi coverings using the \textit{Eisenstein number} which is defined as  $$\displaystyle\sum_{\beta: S \rightarrow \mathbb{P}^{1}} \frac{1}{\text { Aut } \beta}\,,$$ where $\text { Aut } \beta$ denotes the automorphism group of $\beta\,.$ This number becomes $\frac{1}{\text { Aut } \beta}$ for exceptional Belyi coverings.

\begin{thm}
\begin{equation*}
\frac{1}{|\operatorname{Aut} \beta|}=\frac{\left|C_{0}\right|\left|C_{1}\right|\left|C_{\infty}\right|}{(n !)^{2}} \sum_{\chi: \text { irreducible }} \frac{\chi\left(c_{1}\right) \chi\left(c_{2}\right) \chi\left(c_{3}\right)}{\chi(1)} 
\end{equation*}
where $n$ is the degree of $\beta$, $\chi\left(c_{i}\right)$ are the irreducible characters of permutations $c_{i}$ in the ramification scheme and $\left|C_{i}\right|$ are the size of conjugacy classes with representatives $c_{i}$. 
\end{thm}

\section{Elliptic Exceptional Belyi Coverings} 
Elliptic exceptional Belyi coverings are Belyi coverings for the compact Riemann surfaces of genus 1 that are uniquely determined by their ramification schemes.
\subsection{Examples}
\begin{example}
    The following dessin d'enfant is drawn to be like a ``wall-paper". 

\begin{figure}[htp]
\includegraphics[width=50mm]{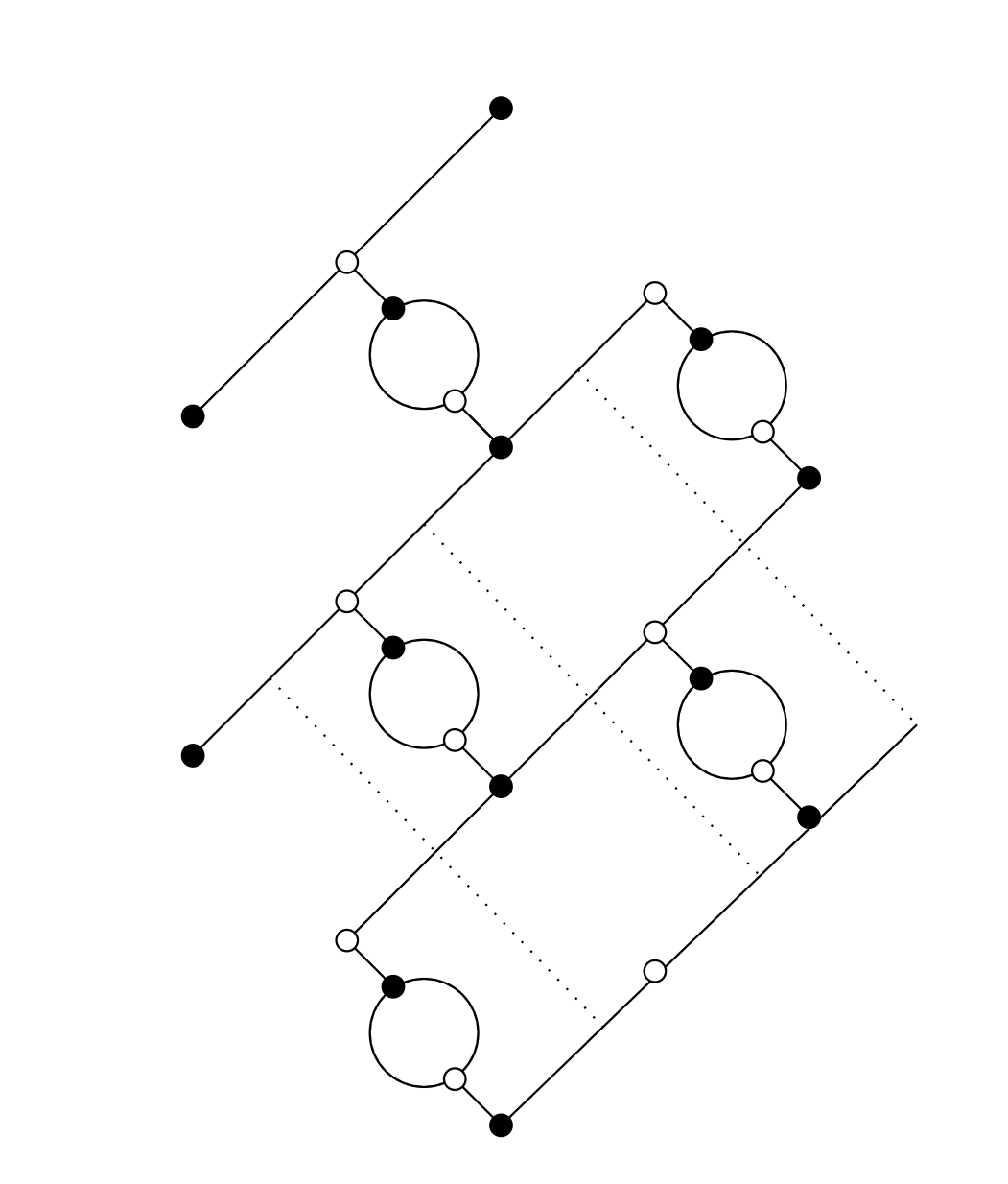}
\end{figure}

We can equivalently think of a compact Riemann surface of genus 1 as a torus. However, trying to draw the respective dessin over this surface is neither an easy nor an illustrative thing to do. Therefore, we consider torus as $\mathbb{C} / \lambda$, where $\lambda$ is a lattice.\\

The respective ramification scheme for this dessin is $[1,5][3,3][3,3]$. This corresponds to the line E6.4 in our table. We also have the information of monodromy permutation:

$$
(1)(2,3,4,5,6)|(4,2,1)(3,5,6)|(1,6,4)(2,3,5) \text {. }
$$
\end{example}

\begin{example}    
     Consider the rational covering corresponding to [3][3][1,1,1] which is given by the map $x \mapsto x^{3}$. It is a star-tree with one black vertex and 3 white vertices. This shows that one black vertex corresponds to 0 with multiplicity 3, the solution of $x^{3}=0$ and each white vertex corresponds to the solutions of $x^{3}=1$, with multiplicity 1. The roots are 1 , $\xi$ and $\xi^{2}$, where $\xi=e^{\frac{2\pi}{3}}$. \\
     
     Now consider the elliptic curve including the square roots of the roots of $x^{3}-1$. Call this curve $E: y^{2}=\left(x^{3}-1\right)$ and consider the projection $E \rightarrow \mathbb{P}^{1}$ by $(x, y) \mapsto x$.

$$
E \rightarrow \mathbb{P}^{1} \rightarrow \mathbb{P}^{1}
$$

The composition in this way yields to a Belyi covering. We call this new covering a ``double-covering" of the rational covering.\\

This corresponds to E6.1 in our table with the ramification scheme $[6][3,3][2,2,2]$. \\
\end{example}

\begin{rem}
 A \textit{field of definition} of a Belyi pair $(S_g, \beta)$, or a dessin denfant, is the smallest number field $K$ such that both the algebraic curve $C$ (corresponding to $S_g$) and $\beta$ can be defined with coefficients in $K$. For example; $\mathbb{Q}(\sqrt{3})$ is the field of definition in {\bf Example \ref{Example 2.1.8}}. In general, existing theories give no upper bound on degree of the field of definition of an exceptional covering when $g=1\,.$  
\end{rem}

\subsection{Maple Calculations}

The Maple code in \cite{Ku} is used to count rational exceptional Belyi coverings. The same code also allows us to determine elliptic exceptional Belyi coverings. However, our knowledge for elliptic coverings is limited in contrast to rational coverings. Our results are limited up to degree 12 and also bounded by ramification schemes, monodromy permutations and Eisenstein numbers. We are not able to list all Belyi pairs for elliptic coverings. Similarly, the respective dessins are not easy to draw as we discussed above. We attach the code and the table for elliptic exceptional Belyi coverings up to degree 12 and we add some descriptions regarding whether they come from modular curves or they are referred as double/triple-coverings. The code shows that there is no elliptic covering of degree 7.  It is an open question whether the number of elliptic exceptional Belyi coverings is finite or infinite. 


\section*{Supplementary material}

```
with(combinat);
with(group);
with(combstruct);
#######################
n:=N: # degree of a covering with three singular fibers
gen:=g: #genus of the covering
print('Genus'=gen,'Degree'=n);
tot:=time():
ABC:=iterstructs(Partition(n+2-2*gen),size=3):##
st:=time():
abc:=nextstruct(ABC): # the number of points in singular fibers
a:=abc[1]:b:=abc[2]:c:=abc[3]:
print('a'=a,'b'=b,'c'=c,'degree'=n);
#n:=a+b+c:       # degree of the covering to ensure that genus=1,
                 # *provided* the covering consists of one component
sol:=0: ## count combinatorail types of covering print('a'=a,'b'=b,'c'=c,degree=n);
unit:=[seq(1,i=1..n)]: ## class of identical permutation
S||n:=permgroup(n, {[[1,2]],[[seq(i,i=1..n)]]}): #symmetric group S_n
P_a := iterstructs(Partition(n),size=a): #ramification indices at fiber a
P_b := iterstructs(Partition(n),size=b): #ramification indices at fiber b
P_c := iterstructs(Partition(n),size=c): #ramification indices at fiber c
while not finished(P_a) do p_a:=nextstruct(P_a):
red_a:=map(x->if x=1 then NULL else x fi,p_a): ## removing cycles of length one
C_a:=SnConjugates(S||n,red_a): # order of the conjugacy class of p_a
while not finished(P_b) do p_b:=nextstruct(P_b):
red_b:=map(x->if x=1 then NULL else x fi,p_b):
C_b:=SnConjugates(S||n,red_b): # order of the conjugacy class of p_b
while not finished(P_c) do p_c:=nextstruct(P_c):
red_c:=map(x->if x=1 then NULL else x fi,p_c):
C_c:=SnConjugates(S||n,red_c): # order of the conjugacy class of p_c
      Rep := iterstructs(Partition(n)): ## Young diagrams = irreps of S_n
      Cov:=0: ## Counting coverings
      while not finished(Rep) do
         rep:=nextstruct(Rep):
Cov:=Cov+Chi(rep,p_a)*Chi(rep,p_b)*Chi(rep,p_c)/Chi(rep,unit):
      od;
         Cov:=Cov*C_a*C_b*C_c/(n!)^2: ## Eisenstein number of coverings

      g:=1/2*(add(x-1, x in p_a)+add(x-1, x in p_b)+add(x-1, x in p_c))-n+1; ## genus, should be = 1.
      if Cov>0 and type(1/Cov,integer)
      then sol:=sol+1:
         print('Ramification indices '= p_a, p_b, p_c);
         print(genus=g,'Eisenstein number of coverings'=Cov):
      fi:
od:
P_c := iterstructs(Partition(n),size=c): #ramification indices at fiber c
od:
P_b := iterstructs(Partition(n),size=b): #ramification indices at fiber b
od:
if sol>0 then print('The number of combinatorial types'=sol, 'time'=time()-st)
else print('No solution','time'=time()-st) fi:
od:
print('Total time'=time()-tot);
```

## Elliptic Exceptional Belyi Coverings

| Label | Ramification scheme | Monodromy permutations | Aut | Description |
|---|---|---|---|---|
| E3.1 | 3 \| 3 \| 3 | (1, 2, 3) \| (1, 2, 3) \| (1, 2, 3) | 3 | $E : y^3 = x(x-1)$; $(x, y) \mapsto x$. |
| E4.1 | 4 \| 4 \| 2, 2 | (4, 1, 2, 3) \| (4, 1, 2, 3) \| (1, 3)(2, 4) | 4 | Double covering of [②]\|[②]\|[①, ①] |
| E4.2 | 4 \| 4 \| 1, 3 | (1, 2, 3, 4) \| (4, 2, 1, 3) \| (1)(2, 3, 4) | 1 | E: |
| E5.1 | 5 \| 5 \| 1, 2, 2 | (5, 1, 2, 3, 4) \| (5, 2, 1, 3, 4) \| (1)(2, 4)(3, 5) | 1 | |
| E5.2 | 5 \| 5 \| 2, 3 | (1, 2, 3, 4, 5) \| (3, 1)(2, 4, 5) \| (3, 5)(1, 2, 4) | 1 | |
| E6.1 | 6 \| 3, 3 \| 2, 2, 2 | (1, 2, 3, 4, 5, 6) \| (3, 5, 1)(2, 4, 6) \| (1, 4)(2, 5)(3, 6) | 6 | Modular; Double covering of [③]\|[3]\|[①, ①, ①] |
| E6.2 | 6 \| 2, 4 \| 2, 2, 2 | (1, 2, 3, 4, 5, 6) \| (1, 4)(2, 3, 5, 6) \| (1, 3)(2, 5)(4, 6) | 2 | Double covering of [③]\|[①, ②]\|[①, 2] |
| E6.3 | 6 \| 1, 5 \| 2, 2, 2 | (1, 2, 3, 4, 5, 6) \| (2)(3, 1, 4, 5, 6) \| (1, 2)(3, 5)(4, 6) | 1 | |
| E6.4 | 1, 5 \| 3, 3 \| 3, 3 | (1)(2, 3, 4, 5, 6) \| (4, 2, 1)(3, 5, 6) \| (1, 6, 4)(2, 3, 5) | 1 | |
| E8.1 | 8 \| 2, 3, 3 \| 2, 2, 2, 2 | (1, 2, 3, 4, 5, 6, 7, 8) \| (1, 5)(2, 4, 7)(3, 6, 8) \| (1, 4)(2, 6)(3, 7)(5, 8) | 2 | Double covering of [④]\|[①, 3]\|[①, ①, 2] |
| E8.2 | 2, 6 \| 4, 4 \| 2, 2, 2, 2 | (1, 2)(3, 4, 5, 6, 7, 8) \| (2, 3, 5, 7)(4, 1, 6, 8) \| (1, 3)(2, 6)(4, 7)(5, 8) | 2 | Double covering of [4]\|[①, ③]\|[①, ①, 2] |
| E8.3 | 4, 4 \| 2, 2, 4 \| 2, 3, 3 | (1, 2, 3, 4)(5, 6, 7, 8) \| (5, 3)(4, 6)(1, 2, 7, 8) \| (1, 7)(2, 4, 5)(3, 8, 6) | 2 | Double covering of [4]\|[①, 3]\|[①, ①, ②] |
| E8.4 | 2, 6 \| 3, 5 \| 2, 2, 2, 2 | (1, 2)(3, 4, 5, 6, 7, 8) \| (2, 3, 6)(4, 1, 5, 7, 8) \| (1, 3)(2, 5)(4, 7)(6, 8) | 1 | |
| E8.5 | 3, 5 \| 3, 5 \| 2, 2, 2, 2 | (1, 2, 3)(4, 5, 6, 7, 8) \| (1, 4, 6)(2, 5, 3, 7, 8) \| (1, 5)(2, 7)(3, 4)(6, 8) | 1 | |
| E8.6 | 1, 7 \| 4, 4 \| 2, 2, 2, 2 | (1)(2, 3, 4, 5, 6, 7, 8) \| (3, 1, 2, 6)(4, 5, 7, 8) \| (1, 2)(3, 5)(4, 7)(6, 8) | 1 | |
| E9.1 | 9 \| 3, 3, 3 \| 1, 2, 2, 2, 2 | (1, 2, 3, 4, 5, 6, 7, 8, 9) \| (4, 8, 2)(6, 5, 1)(3, 7, 9) \| (5)(1, 4)(2, 7)(3, 8)(6, 9) | 1 | |
| E9.2 | 4, 5 \| 3, 3, 3 \| 2, 2, 2, 3 | (1, 2, 3, 4)(5, 6, 7, 8, 9) \| (5, 7, 2)(3, 6, 1)(4, 8, 9) \| (1, 5)(2, 6)(7, 9)(3, 4, 8) | 1 | |
| E9.3 | 3, 3 \| 3, 3, 3 \| 2, 2, 5 | (1, 2, 3)(4, 5, 6)(7, 8, 9) \| (4, 2, 7)(1, 5, 8)(3, 6, 9) \| (1, 7)(2, 6)(3, 8, 4, 9, 5) | 1 | |
| E10.1 | 5, 5 \| 3, 3, 4 \| 2, 2, 2, 2, 2 | (1, 2, 3, 4, 5)(6, 7, 8, 9, 10) \| (3, 7, 1)(2, 6, 8)(4, 5, 9, 10) \| (1, 6)(2, 7)(3, 5)(4, 9)(8, 10) | 2 | Double covering of [5]\|[①, 3]\|[①, ①, ①, 2] |
| E10.2 | 5, 5 \| 1, 4, 5 \| 2, 2, 2, 2, 2 | (1, 2, 3, 4, 5)(6, 7, 8, 9, 10) \| (6, 9, 5, 2)(3, 4, 7, 1, 8)(10) \| (1, 6)(2, 4)(3, 7)(5, 8)(9, 10) | 1 | |
| E10.3 | 3, 7 \| 2, 4, 4 \| 2, 2, 2, 2, 2 | (1, 2, 3)(4, 5, 6, 7, 8, 9, 10) \| (2, 5)(3, 4, 7, 9)(6, 1, 8, 10) \| (1, 5)(2, 4)(3, 8)(6, 9)(7, 10) | 1 | |
| E10.4 | 3, 7 \| 3, 3, 4 \| 2, 2, 2, 2, 2 | (1, 2, 3)(4, 5, 6, 7, 8, 9, 10) \| (6, 9, 3)(4, 7, 2)(5, 1, 8, 10) \| (1, 4)(2, 6)(3, 8)(5, 9)(7, 10) | 1 | |
| E10.5 | 10 \| 1, 3, 3, 3 \| 2, 2, 2, 2, 2 | (1, 2, 3, 4, 5, 6, 7, 8, 9, 10) \| (4, 7, 2)(3, 6, 8)(5, 1, 9)(10) \| (1, 4)(2, 6)(3, 7)(5, 8)(9, 10) | 1 | |
| E12.1 | 1, 11 \| 3, 3, 3, 3 \| 2, 2, 2, 2, 2, 2 | (1, 2, 3, 4, 5, 6, 7, 8, 9, 10, 11)(12) \| (4, 7, 2)(3, 6, 8)(5, 1, 9)(10, 11, 12) \| (1, 4)(2, 6)(3, 7)(5, 8)(9, 11)(10, 12) | 1 | Modular; $y^2 + y = x^3 - x^2 - 10x - 20$ |
| E12.2 | 2, 10 \| 3, 3, 3, 3 \| 2, 2, 2, 2, 2, 2 | (1, 2)(3, 4, 5, 6, 7, 8, 9, 10, 11, 12) \| (2, 3, 9)(5, 7, 11)(8, 4, 1)(6, 10, 12) \| (1, 3)(2, 8)(4, 7)(5, 10)(6, 11)(9, 12) | 2 | Modular; 2-covering of [①, ⑤]\|[3, 3]\|[①, 2, 2] |
| E12.3 | 3, 9 \| 3, 3, 3, 3 \| 2, 2, 2, 2, 2, 2 | (1, 2, 3)(4, 5, 6, 7, 8, 9, 10, 11, 12) \| (4, 8, 2)(5, 1, 10)(3, 7, 11)(6, 9, 12) \| (1, 4)(2, 7)(3, 10)(5, 9)(6, 11)(8, 12) | 3 | Modular; 3-covering of [①, ③]\|[①, 3]\|[2, 2], see E3.1 |
| E12.4 | 4, 8 \| 3, 3, 3, 3 \| 2, 2, 2, 2, 2, 2 | (1, 2, 3, 4)(5, 6, 7, 8, 9, 10, 11, 12) \| (5, 8, 2)(3, 7, 10)(6, 1, 11)(4, 9, 12) \| (1, 5)(2, 7)(3, 9)(4, 11)(6, 10)(8, 12) | 4 | Modular; 4-covering of [①, ②]\|[①, 3]\|[①, 2], see E4.1 |
| E12.5 | 5, 7 \| 2, 2, 4, 4 \| 2, 2, 2, 2, 2, 2 | (1, 2, 3, 4, 5)(6, 7, 8, 9, 10, 11, 12) \| (8, 2)(6, 10)(1, 3, 7, 9)(4, 5, 11, 12) \| (1, 8)(2, 7)(3, 9)(4, 11)(6, 9)(10, 12) | 1 | |
| E12.6 | 6, 6 \| 1, 3, 4, 4 \| 2, 2, 2, 2, 2, 2 | (1, 2, 3, 4, 5, 6)(7, 8, 9, 10, 11, 12) \| (2)(5, 7, 9)(11, 6, 8, 4)(3, 1, 10, 12) \| (1, 2)(3, 11)(4, 7)(5, 8)(6, 10)(9, 12) | 1 | |
| E12.7 | 2, 3, 7 \| 3, 3, 3, 3 \| 2, 2, 2, 2, 2, 2 | (1, 2)(3, 4, 5)(6, 7, 8, 9, 10, 11, 12) \| (5, 2, 8, 11)(4, 6, 9, 1)(7, 3, 10, 12) \| (1, 8)(2, 4)(3, 6)(5, 10)(7, 11)(9, 12) | 1 | |
| E12.8 | 3, 3, 6 \| 3, 4, 5 \| 2, 2, 2, 2, 2, 2 | (1, 2, 3)(4, 5, 6)(7, 8, 9, 10, 11, 12) \| (3, 4, 11)(7, 9, 5, 2)(6, 8, 1, 10, 12) \| (1, 7)(2, 4)(3, 10)(5, 8)(6, 11)(9, 12) | 1 | |



\section{Acknowledgements}
I thank my supervisor Professor Alexander Klyachko for his valuable support and guidance.

\bibliographystyle{alpha}
\bibliography{mybib}

\end{document}